\documentclass[a4paper,fontsize=9pt]{article}

\usepackage[english]{babel}
\usepackage[utf8]{inputenc}
\usepackage{amsmath}
\usepackage{graphicx}
\usepackage[colorinlistoftodos]{todonotes}
\usepackage{psfrag,epsf}
\usepackage{enumerate}
\usepackage[numbers]{natbib}
\usepackage{url} 
\usepackage{amsthm,mathrsfs,amsfonts}	
\usepackage{enumitem}
\usepackage{bbm}

\numberwithin{equation}{section}
\newtheorem{theorem}{Theorem}[section]
\newtheorem{definition}{Definition}[section]

\newtheorem{remark}{Remark}

\newtheorem{proposition}{Proposition}

\newcommand{\R}{\mbox{$\mathbb{R}$}}

\newcommand{\ind}{\mathbb{I}}

\newcommand{\C}{\mathcal{C}}

\usepackage{graphicx} 
\usepackage{subfig}
\usepackage{float}		
\usepackage{color}
\usepackage[normalem]{ulem}
\title{\bf Level set and density estimation on manifolds}

\author{Alejandro Cholaquidis \thanks{acholaquidis@cmat.edu.uy}\hspace{.2cm}\\
	Centro de Matem\'atica, Facultad de Ciencias, \\
	Universidad de la Rep\'ublica, Uruguay.\\
	Ricardo Fraiman  \\
	Centro de Matem\'atica, Facultad de Ciencias,\\
	Universidad de la Rep\'ublica, Uruguay.\\
	and \\ 
	Leonardo Moreno \\ 
	Departamento de M\'etodos Cuantitativos, FCEA, \\
	Universidad de la Rep\'ublica, Uruguay.}
\begin{document}
	\maketitle

\begin{abstract}
	We tackle the problem of the estimation of the level sets $L_f(\lambda)$ of the density $f$ of a random vector $X$ supported on a smooth manifold $M\subset \mathbb{R}^d$, from an iid sample of $X$. To do that  we introduce a kernel-based estimator $\hat{f}_{n,h}$, which is a slightly modified version of the one proposed in \cite{snrc:14}, and proves its a.s. uniform convergence to $f$. Then, we propose two estimators of $L_f(\lambda)$, the first one is a plug-in: $L_{\hat{f}_{n,h}}(\lambda)$, which is proven to be a.s. consistent in Hausdorff distance and distance in measure, if $L_f(\lambda)$ does not meet the boundary of $M$. While the second one assumes that $L_f(\lambda)$ is $r$-convex, and is estimated by means of the $r$-convex hull of $L_{\hat{f}_{n,h}}(\lambda)$.
    The performance of our proposal is illustrated through some simulated examples. In a real data example we analyze the intensity and direction of strong and moderate winds.
\end{abstract}


\section{Introduction}

The statistics of functional data had an exponential growth in the last decades, being particularly important the impulse given by the early works \cite{ram02} and \cite{ram05}. Very different problems have been addressed in recent years, see for instance the reviews by \cite{cue14}, \cite{go:16} \cite{an:19}. Initially the results were concentrated on problems where the data took values in functional spaces, such as $L^2, L^1$ and $L^{\infty}$, and more generally on Hilbert or Banach spaces. Recently a lot of interest appears for data in more general spaces (where there is no linear structure on them), like metric spaces. This is motivated by important practical applications, see for instance subsection 8.1. Some relevant examples includes the case of high dimensional data (see for instance \cite{has89}, \cite{del01}), directional data (see \cite{mar00}, \cite{mar72}), cone and cylindric data (see example on subsection 8.1), and random graph data (see  \cite{frai}), among others.

In what follows we will consider the case where the data take values on a Riemannian manifold.

 Starting from the pioneer 1945 work of Rao (see \cite{rao45}), the statistical theory for data valued on a Riemannian manifold has received a lot of interest  because of its important applications.  In particular, these techniques may allow to avoid the curse of dimensionality when trying to analyze data in a high dimensional ambient space. 
Indeed, as mentioned in \cite{hen:07}:  ``\textit{Data belonging to some $m$-dimensional compact submanifold $M$ of Euclidean space $R^s$ appear in many areas of natural science. Directional statistics, image analysis, vector cardiography in medicine, orientational statistics, plate tectonics, astronomy and shape analysis comprise a (by no means exhaustive) list of examples}''. These techniques are also applied in medical imaging applications: as it mentioned in \cite{pennec} ``\textit{Examples of manifolds we routinely use in medical imaging applications are 3D rotations, 3D rigid transformations, frames (a 3D point and an orthonormal trihedron), semi- or non-oriented frames [...]  , positive definite symmetric matrices coming from diffusion tensor imaging}''.

The estimation of level sets $L_f(\lambda) = \{x : f (x) \geq  \lambda\}$, where
$f$ is an unknown density function on $\mathbb{R}^d$ and $\lambda> 0$ is a given constant, has been considered by many authors; see, for instance, \cite{harti87}, \cite{polo95}, \cite{cuefrai97}, \cite{mol98}
\cite{tsy97}, \cite{wal97} for consistency results and rates of convergence, while the asymptotic distribution was derived in  \cite{wass17}. Some relevant applications include mode estimation  \cite{mull91}, \cite{polo95}, clustering (\cite{cuefrai00}, \cite{cuefrai01}) or detection of abnormal behaviour in a system (\cite{dw80}, \cite{bca01}, \cite{b03}). However, this problem is less developed when the underlying density has its support on a Riemannian manifold. In the following we address the problem of level set estimation in this setup.

More precisely, given a $d'$-dimensional Riemannian manifold $M\subset \mathbb{R}^d$, where $M$ is unknown but $d'\leq d$ is assumed to be known, the aim is to estimate $L_f(\lambda)$ of the density $f$ of a random vector $X$ with support $M$ from an iid sample $X_1,\dots,X_n$ of $X$. In practice and for large sample sizes, if the intrinsic dimension of $M$ is unknown it can be first  estimated by means of, for instance, the classical estimator proposed in  \cite{lb05}.

Our first proposed estimator is just the plug-in estimator $L_{\hat{f}_{n,h}}(\lambda)$, where $\hat{f}_{n,h}$ is a kernel-based estimator of $f$ with bandwidth $h=h_n\to 0$, which is a slightly modified version of the one proposed in \cite{snrc:14}. The almost sure (a.s.) consistency of $L_{\hat{f}_{n,h}}(\lambda)$ requires to prove the a.s. uniform convergence of the sequence of estimators $\hat{f}_{n,h}$ to $f$, this is done in Section \ref{densest}. Regarding density estimation on manifolds, in \cite{pelletier2006}  $L^2$-consistency is obtained for a kernel-based density estimator (w.r.t. the $L^2$ norm in $M$, see Theorem 3.1).   There are two main drawbacks to use that result to our setup, first in \cite{pelletier2006} it is assumed that the Riemannian structure in $M$ is known, and second, $L^2$ consistency does not imply convergence of level sets.  For the same estimator, \cite{hr09} obtain the limit distribution and the a.s. consistency for the uniform metric. In \cite{jiang2017b} a different kernel-based estimator is proposed (similar to the one we will propose), and convergence in probability is obtained for the uniform metric. It does not assume that the underlying Riemannian structure is known. All the aforementioned results are for manifolds without boundary. 	Other references that tackle the density estimation problem for manifolds without boundary are \cite{ozgr09} and \cite{wass19}. For manifolds with boundary point-wise $L^2$ consistency is obtained in \cite{bs17}. 
Lastly, in section \ref{rconv} we tackle the level set estimation problem but imposing a well known-shape restriction called $r$-convexity. In this case the proposed estimator is the $r$-convex hull of $L_{\hat{f}_{n,h}}$.

\section{Roadmap}

 In Section \ref{notacion} we introduce the basic notation and the geometric framework used throughout the manuscript. Section \ref{densest} is devoted to prove the a.s. uniform convergence of $\hat{f}_{n,h}$ to $f$. In Section \ref{levest} we prove that if $L_f(\lambda)$ does not meet the boundary of $M$,  $L_{\hat{f}_{n,h}}(\lambda)$ converges a.s. in Hausdorff distance as well as in distance in measure to $L_f(\lambda)$. We also prove the convergence of its boundaries, i.e: $\partial L_{\hat{f}_{n,h}}(\lambda)$ converges in Hausdorff distance to $\partial L_f(\lambda)$. If $L_f(\lambda)$ meets the boundary of $M$, we prove that $L_{\hat{f}_{n,h}}(\lambda)$ converges a.s. in Hausdorff distance to $L_f(\lambda)$.  Consistency in the Hausdorff metric of level sets under $r$--convexity is shown in Section 6. In Section 7  we provide some simulation results, while in Section \ref{winds} we consider an important application to describe the wind behavior in Uruguay. All proofs are given in the appendix.

\section{Notation and geometric framework}\label{notacion}

If $B\subset\mathbb{R}^d$ is a Borel set, then we denote by $|B|$ its Lebesgue measure and by $\overline{B}$ 
its closure.   Given a set $A$ on a topological space, the interior of $A$ with respect to the underling topology is denoted by $\mathring{A}$.
The $k$-dimensional closed ball of radius $\varepsilon$ centered at $x$ will be denoted by
$\mathcal{B}_k(x,\varepsilon)\subset \mathbb{R}^d$ (when $k=d$ the index will be omitted), and 
its Lebesgue measure is  denoted by $\sigma_k=|\mathcal{B}_k(x,1)|$. The Euclidean inner product in $\mathbb{R}^d$ is denoted by $\langle \cdot ,\cdot\rangle$, while the corresponding norm in $\mathbb{R}^d$ is denoted by $\|\cdot\|$.


From now on, we assume that $M\subset \mathbb{R}^d$ is a compact  $d'$-dimensional manifold of class 
$\mathcal{C}^2$ (also called a $d'$-regular surface of class $\C^2$). We consider the Riemannian metric on $M$ inherited
from $\mathbb{R}^d$. If $x\in M$, $T_xM$ denotes the tangent space at $x$, while $\rho(x,y)$ denotes the geodesic distance between $x$ and $y$.  Given a set $A\subset M$, we denote $B_\rho(A,r)=\{x\in M: \rho(x,A)<r\}$.  For $f:M\to \mathbb{R}$, we denote by $\nabla f(x)$ the gradient of $f$ at $x\in M$.
When $M$ is orientable, it has a unique associated volume form $\omega$ such that $\omega(e_1,\ldots, e_{d'})=1$ 
for all oriented orthonormal
bases $e_1,\ldots, e_{d'}$ of $T_xM$. If $g:M\rightarrow \mathbb{R}$ is a density function, then
we can define a new measure
$\mu(B)=\int_B gd\omega$, where $B\subset M$ is a Borel set.  In what follows we assume that $M$ is orientable.
Given a point $x\in M$, $b_x$ is the geodesic distance from $x$ to the boundary $\partial M$ of $M$, or is $\infty$ if $\partial M=\emptyset$.  

Recall that given two non-empty compact sets $A,C\subset \mathbb{R}^d$, the Hausdorff distance between $A$ and $C$ is defined as 
\begin{equation}
d_H(A,C)=\max\Big\{\max_{a\in A}\rho(a,C), \ \max_{c\in C}\rho(c,A)\Big\}, \text{ where }
\rho(a,C)=\inf_{c\in C} \rho(a,c).
\end{equation}

Given two Borel sets $A,B\subset M$, the distance in measure between them is $d_\mu(A,B)=\mu(A\setminus B)+\mu(B\setminus A)$.


\section{Density estimation} \label{densest}

The aim of this section is to prove that $\hat{f}_{h,n}$, a modified version of the kernel-based density estimator, denoted by $f_{h,n}$, proposed in \cite{bs17},  converges uniformly to the density $f$ when the manifold has a $\C^2$ boundary.  This auxiliary result, besides the interest in itself, will be used to prove our main results regarding level set estimation in the next section. Let us recall the definition of $f_{h,n}$.  We assume that $K$ is the Gaussian kernel, (however, it can be replaced by any sub Gaussian distribution, see  Remark \ref{rem1} below); that is, $K(\|x\|)=\pi^{-d'/2}\exp(-\|x\|^2)$. Let $h=h_n\to 0$; then,  
\begin{equation}\label{kerndens}
f_{h,n}(x)= \frac{1}{nm_0(x)h^{d'}}\sum_{i=1}^n K\Big(\frac{\|x-X_i\|}{h}\Big)\ \, \textrm{where} \ \, m_0(x)=\pi^{-1/2}\int_{-\infty}^{b_x/h} \exp(-z^2)dz,
\end{equation}
where $b_x$ is the distance to $\partial M$ of the point $x \in M$.  
An important assumption in \cite{bs17}, requires that
the manifold $M$  be ``uniformly tangible", which, roughly speaking, allows to define the projection onto the boundary of points close enough to it. More precisely:

\begin{definition} A $d'$-dimensional Riemannian manifold $M\subset \mathbb{R}^d$   is said to be uniformly tangible if
\begin{enumerate}
	\item[1.] There exists $\delta>0$ and $c>0$ such that $R(x,y):=\|x-y\|/\rho(x,y)>c$ for all $x,y\in M$, such that $\|x-y\|<\delta$.
	\item[2.] There exists $r>0$ such that   $N_r=\partial M\times[0,r)$ is mapped diffeomorphically onto its image via the exponential map $(x,t)\to \exp_x(tv_x)$ where $x\in \partial M$,  $t\in[0,r)$, and $v_x$ is the inward pointing unit normal vector to
the boundary.
	\item[3.] Denote by $\text{inj}(x)$ the injectivity radius of a point x (i.e the maximum radius for which $B(0,r)\subset T_xM$ is mapped diffeomorphically into $M$ by the exponential map).  Then $\inf\{\text{inj}(x): x\notin \exp(N_r)\}>0$.
\end{enumerate}
\end{definition}

The following proposition states that this condition holds when $M$ is a compact $\C^2$ manifold whose boundary (in case there exists) is also a $\C^2$ manifold. Under this condition there exists a radius $r_M>0$ such that for any point $x$ within a uniform geodesic distance $r_M$ to the boundary, there exists a unique closest point on $\partial M$. Then it can be defined $\eta_x$ the unit vector pointing in the direction of the unique closest boundary point. For points farther away than $r_M$ from the boundary, $\eta_x$ can be chosen arbitrarily.

\begin{proposition} \label{unitang} Let $M\subset \mathbb{R}^d$ be a $d'$-dimensional $\C^2$ compact Riemannian manifold, whose boundary, in case there exists, is a $\C^2$ manifold. Then, $M$ is uniformly tangible.
	
\end{proposition}

Equation (5) in \cite{bs17} states that,  if $M$ is uniformly tangible, and $h$ is small enough the bias of $f_{h,n}(x)$ is
\begin{equation}\label{bias}
E(f_{h,n}(x))-f(x)= hm_1(x)\langle \eta_x , \nabla f(x)\rangle+\mathcal{O}_x(h^2) \ \  \text{ where } m_1(x)=\frac{1}{2\sqrt{\pi}}\exp(-b_x^2/h^2).
\end{equation}

 Observe that if $\lambda$ is such that  $L_f(\lambda)\cap \partial M\neq \emptyset$, then if $x\in L_f$, $m_0(x)\to 1$. This suggest to use the following estimator:

\begin{equation} \label{nuestroest}
\hat{f}_{h,n}(x)= \frac{1}{n h^{d'}}\sum_{i=1}^n K\Big(\frac{\|x-X_i\|}{h}\Big),
\end{equation}
  which does not depend on $m_0(x)$ or $b_x$.  

In order to get the a.s. convergence of $L_{\hat{f}_{h,n}}(\lambda)$ to $L_f(\lambda)$ we need to prove the following auxiliary result that states the a.s. uniform convergence of $\hat{f}_{h,n}$ to $f$.

\begin{theorem} \label{convunif2} Under the hypotheses of Proposition \ref{unitang}. Let $X$ be a random vector with support $M$ whose density $f$ is assumed to be $\mathcal{C}^4$. Let $h\to 0$ such that $nh^{d'+3}/\log(n)\to \infty$ when $n \to +\infty$; then, 
	$$\sup_{x\in M_n} |\hat{f}_{h,n}(x)-f(x)|=o(h/c_n) \quad a.s.$$
	for any sequence of closed subsets $M_n\subset M$ such that $c_n/h\to \infty$, where $c_n=\inf_{x\in M_n}\rho(x,\partial M)>0$.
\end{theorem}

The following result is more restrictive, but a better rate of convergence is obtained,  (it holds in particular when $\partial M=\emptyset$) since the theorem does not allows the compact set $M_0$ to depend on $n$. This is proven in the same manner as Theorem \ref{convunif2},

\begin{theorem} \label{convunif} Under the hypotheses of Theorem \ref{convunif2}. Let $h\to 0$  and $\beta_n\to \infty$ such that $\beta_n h^2\to 0$, $nh^{d'}/(\beta_n^2\log(n))\to \infty$; then, 
	$$\beta_n\sup_{x\in M_0} |\hat{f}_{h,n}(x)-f(x)|\to 0\quad a.s.$$
	for any closed subset $M_0\subset M$ such that $\inf_{x\in M_0}\rho(x,\partial M)>0$.
\end{theorem}
 
\begin{remark}
	As a consequence, taking $\beta_n = n^{\alpha}$ and $h= n^{-\gamma}$ we derive that we can reach $\beta_n = n^{\alpha}$ for any $\alpha < 2/(d'+4)$.
\end{remark}
\begin{remark} \label{rem1} The explicit expression for $m_0$ and $m_1$ given in \eqref{kerndens} and \eqref{bias} respectively, strongly rely on the gaussianity of the kernel. 	If  this assumption is removed, the expressions are much more involved, see \cite{bs17}. Observe that the proposed estimator \eqref{nuestroest} does not require the computation of $m_0$ or $m_1$, but the proofs of Theorems  \ref{convunif2} and \ref{convunif} uses them as auxiliary tools. However, the same rate of convergence is obtained in Theorems \ref{convunif2} and \ref{convunif} if the gaussian kernel is replaced by any Lipchitz kernel fulfilling $K(\|z\|)\leq C_1\exp(-C_2\|z\|^2)$, with a subgaussian density function, i.e. fulfilling $\mathbb P(\|X\| > t ) \leq C_1\exp(-C_2\|t\|^2)$, %
 for all $t$ and, for some positive constants $C_1,C_2$, where $X$ is a $d'$-dimensional vector with density  $K(\|z\|)$.
\end{remark}

\section{Level set estimation}\label{levest}

Level set estimation is an important problem with many applications in statistics, such as in hierarchical clustering, binary classification, outliers detection, functional neuroimaging,  and bioinformatics among many others. References are given in the introduction. In our setup we consider this problem when the distribution of the data is supported on a smooth Riemannian manifold $M \subset \R^d$ which as mentioned in the introduction covers important applications that includes directional data, cone and cilindrical data,  high dimensional data and random graph data. In what follows we state our main asymptotic results.  
   
Once we prove the a.s. uniform convergence of the estimator $\hat{f}_{n,h}$ to $f$ we are ready to state our main results.  The first one (Theorem \ref{th:3}) tackle the case in which $L_f(\lambda)$ does not meet the boundary of the manifold. For this case we obtained not only the a.s. convergence in Hausdorff distance of the level set estimator $L_{\hat{f}_{n,h}}(\lambda)$, but also the convergence of its boundary, as well as the convergence in measure. As  usual in level set estimation, we require that the boundary of $\partial L_f(\lambda)$, which is $\{f=\lambda\}$ because we assume that $f$ is a continuous function, does not contain a plateau at level $\lambda$, see the discussion on condition f1) in Theorem 1 in \cite{cgmrc:06}. The proof is based on Theorem \ref{convunif}.

\begin{theorem} \label{th:3}  Let $M$ and $f$ as in Theorem \ref{convunif2}. Assume that the level $\lambda>0$ fulfills that for all $x$ such that $f(x)=\lambda$, there exists $a_n,b_n\to x$ such that $f(a_n)>\lambda $ and $f(b_n)<\lambda$ and   $\partial L_f(\lambda)\cap \partial M=\emptyset$. Then, with probability one,
	\begin{itemize}
		\item[1.] 	$d_H(\partial L_{\hat{f}_{n,h}}(\lambda),\partial L_{f}(\lambda))\to 0$;
		\item[2.]  $d_H(L_{\hat{f}_{n,h}}(\lambda),L_f(\lambda))\to 0$;
		\item[3.] If, moreover, $\nabla_x f\neq 0$ for all $x$ such that $f(x)=\lambda$, $d_\mu(L_{\hat{f}_{n,h}}(\lambda),L_f(\lambda))\to 0$.
	\end{itemize}
\end{theorem}
 
If the underlying level set $L_f(\lambda)$ meets the boundary of the manifold, then we have the following result, whose proof is based on Theorem \ref{convunif2}.
 
\begin{theorem} \label{convhaus} Let $M$ and $f$ be as in  Theorem \ref{convunif2}. Assume that the level $\lambda>0$ fulfills that for all $x$ with $f(x)=\lambda$, there exists $a_j\to x$, $a_j\in \mathring{M}$, such that $f(a_j)>\lambda$ for all $j$. Then,
	$$d_H\big(L_{\hat{f}_{n,h}}(\lambda),L_{f}(\lambda)\big)\to 0, \quad a.s.$$
\end{theorem}

\section{Manifold level set estimation under r-convexity}\label{rconv}
In what follows we consider the level set estimation problem when we assume that the level set is $r$-convex set in the manifold.   

In Euclidean space, a set $A$ is said to be $r$-convex (for some $r>0$) if $A=C_r(A)$, where $C_r(A)$ is the $r$-convex hull of $A$; that is, the intersection of the complements of all open balls of radius $r$ that does not meet $A$. This is a natural generalization of convexity (the half spaces are replaced by balls), and it has been widely studied in set estimation  literature (see, for instance, \cite{wal97,wal99} \cite{rod07} and \cite{pat08}). Additionally, as is pointed out in \cite{rod07}, this concept ``is closely related to the notion of alpha-shapes that arises in the literature of computational geometry"; see \cite{edel94}. Departing from the idea of $r$-convexity, several generalizations have been given (see, for instance, \cite{chola:14}). 
If the underlying space is not Euclidean space but is rather any Riemannian manifold $M$ endowed with the geodesic distance $\rho$, then  the natural generalization  is to replace the Euclidean balls with geodesic balls. According to this idea, given $r>0$, we will say that a set $A\subset M$ is $r$-convex if it is equal to its $r$-convex hull in $M$, that is, the intersection of the complement of all open geodesic balls of radius $r$ that does not meet $A$.

\begin{theorem}\label{th:5} Under the hypotheses of Theorem \ref{convhaus}, assume also that the level sets $L_f(\lambda)$ is $r$-convex and $L_{\hat{f}_{n,h}}(\lambda)$  is $r$-convex a.s., for some $r>0$.
	Then,
	
	$$d_H\big(C_r(\{X_i:\hat{f}_{h,n}(X_i)>\lambda\}),C_r(\{X_i:f(X_i)>\lambda\})\big)\to 0, \quad a.s.$$
	and	
	$$d_H\big(C_r(\{X_i:\hat{f}_{h,n}(X_i)>\lambda\}),L_f(\lambda)\big)\to 0\quad a.s.$$
\end{theorem}

\section{Simulation results}

To assess the performance of our proposal, we will perform a simulation example with two scenarios. In the first one, we consider a distribution on the positive cone of covariance $2 \times 2$-matrices, which is a three dimensional manifold when endowed with the Riemannian structure given below. In the second one, we compare our density estimator with the one proposed in \cite{hall}, which is specially designed for spherical data. Also, as illustrative examples, we consider the torus with the metric inherit from $R^3$ and the two dimensional half-sphere in $\mathbb{R}^3$. In the first case we consider two distributions: the first is unimodal and the second is a mixture of distributions. In the last case we considered a Von-Mises distribution.

\subsection{Positive-definite matrices}

Let us denote by $(\mathbb{P}_d,g)$ the set of positive-definite $d\times d$-covariance matrices. Given two matrices  $A,B \in \mathbb{P}_d$, the geodesic curve joining $A$  and  $B$ is  

\begin{displaymath}
\gamma(s) = A ^{1/2}  (A ^{-1/2} B A ^{-1/2}  ) ^ s A ^{1/2} \quad \text{for all }s\in [0,1].
\end{displaymath}
The geodesic distance is given by $d_g(A,B)=\Vert \ln  (A ^{-1/2} B A ^{-1/2} ) \Vert$, where   $\Vert \cdot \Vert$  is the Hilbert--Schmidt norm.
%
%

\begin{figure}[h]
	\begin{center}
		\includegraphics[scale=0.4]{ 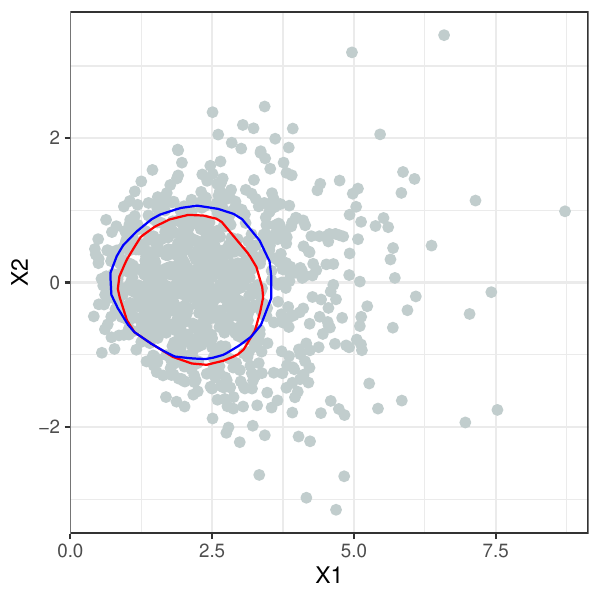}
		 		\includegraphics[scale=0.4]{ 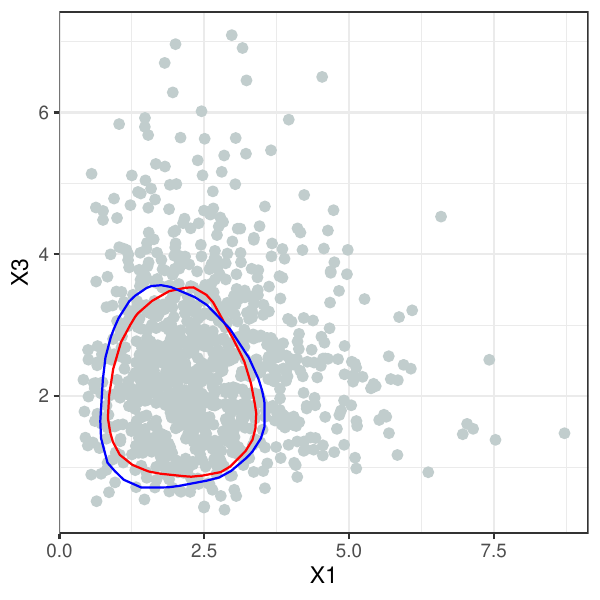}
		 				\includegraphics[scale=0.4]{ 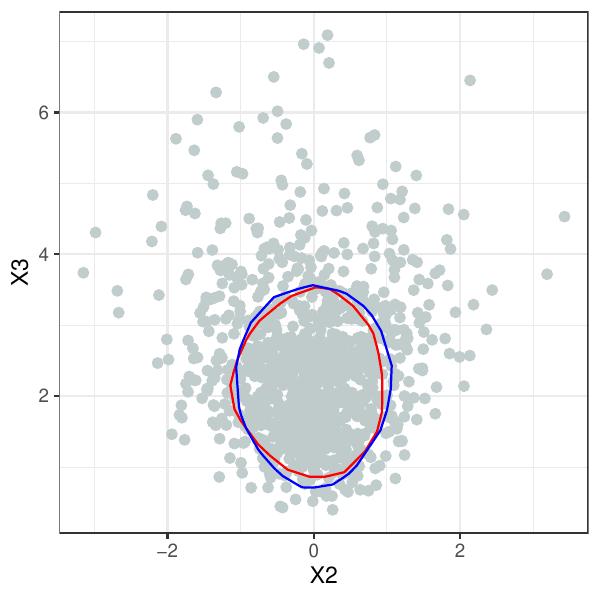}
		\caption{Projections of a sample of size 1000 drawn from a Wishart distribution with $m=10$ and $\Sigma=(1/4)I$, together with the convex hull of the $\lambda$ level set (in blue) and the convex hull of the level set estimator (in red), for $\lambda=0.03$ and $h=0.55$.}
		\label{matrices}
	\end{center}
\end{figure}

We consider, for $d=2$, the Wishart distribution  $\mathcal{W}_2(\Sigma, m)$ on $\mathbb{P}_2$ with parameters $m=10$ and $\Sigma=(1/4)I_2$. An easy way to obtain a matrix $S$ with this distribution is to define $S =  X_1 X'_1 + \cdots + X_m X'_m $, where $X_1,\dots,X_m$ is an iid random sample of a multivariate Gaussian distribution with mean $0$ and covariance matrix $\Sigma$. \\
As is well-known, $(\mathbb{P}_2,g)$ can be represented as a cone in $\mathbb{R}^3$. In Figure \ref{matrices}, we show the projections of a sample of size $1000$, drawn from a Wishart distribution with $m=10$ and $\Sigma=(1/4)I$, together with the convex hull of the $\lambda$ level set $L_{\mathcal{W}}(\lambda)$  (in blue) and the convex hull of the level set estimator $L_{\hat{\mathcal{W}}_{n,h}}(\lambda)$ (in red) for $\lambda=0.03$ and $h=0.55$.  The Hausdorff distance between the level sets in this case is $\mathbb{R}^3$ is $0.325$.
 In Table \ref{power_norm}, we report the mean, median, and standard deviation, over $100$ replications of the Hausdorff distance ($d_H$) between both sets for different sample sizes $n\in \{1000,5000,10000,20000\}$.  The parameter $h$ is chosen following the proposal in Appendix B of \cite{bs17}. 

  \begin{table}[ht]
\begin{center}
\scalebox{1}{
\begin{tabular}{ r| c|c}
	$n$   & $h$  & $d_H$ \\ \hline
	1000  & 0.55  & 0.3388 [0.3286] (0.0616)  \\
	5000  & 0.4   & 0.2326 [0.2313] (0.0266) \\
	10000 & 0.35  & 0.2114 [0.2110] (0.0267) \\
	20000 & 0.3  & 0.1909 [0.1872] (0.0235)
\end{tabular}}
\caption{  Mean over 100 replications for the Hausdorff distance between the true level set  $L_{\mathcal{W}}(\lambda)$ and the estimator $L_{\hat{\mathcal{W}}_{n,h}}(\lambda)$ for $\lambda=0.03$. Between $[]$ and $()$ the median and standard deviation respectively.}
\label{power_norm}
\end{center}
\end{table}

\subsection{The torus}

\begin{figure}[h]
	\begin{center}
		\includegraphics[scale=.27]{ 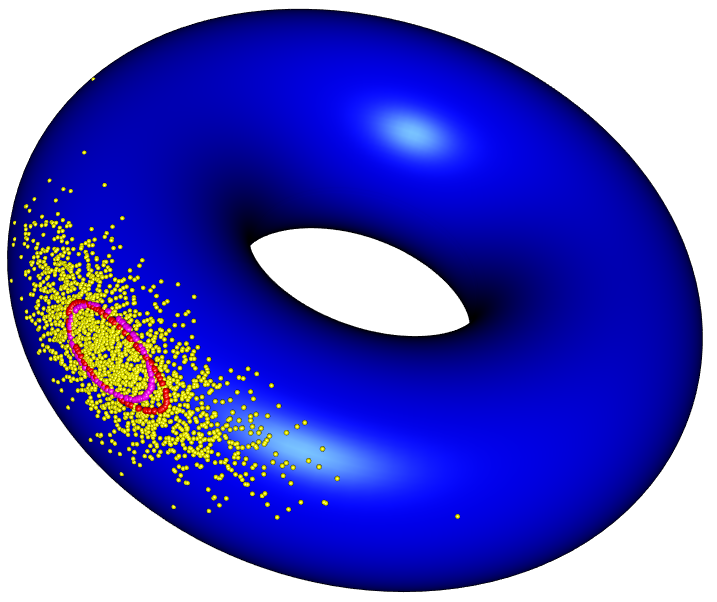}
			\includegraphics[scale=.25]{ 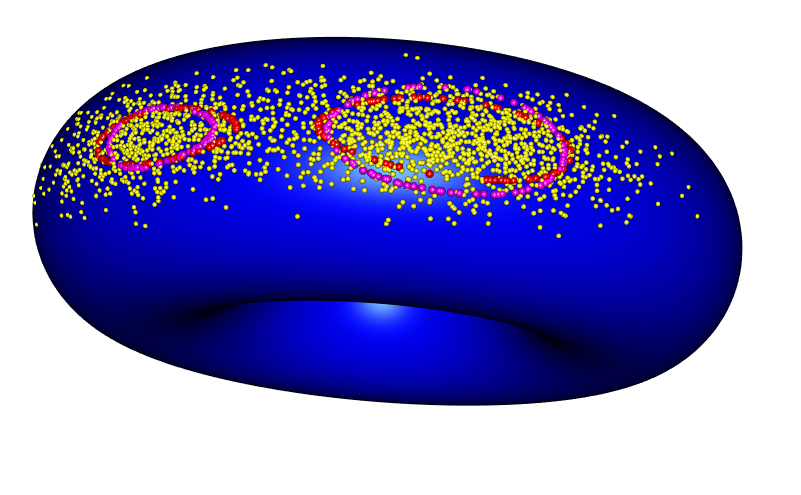}		
		\caption{Left-hand panel:  a sample of size 2000 from a $\mathcal{MVM}_1(\mu_1,\kappa_1, \Delta_1)$ distribution with $\mu_1$, $\kappa_1$ and $\Delta_1$ given in \eqref{law1}. Right-hand panel:  a sample of size 2000 from the  mixture law given in \eqref{mixture}. In both cases, the data are shown in yellow, whereas the boundary of the true level sets is shown (in red) together with the estimated boundary (in magenta).  }
		\label{toro}
	\end{center}
\end{figure}

In the torus $\mathbb{T}^2= S^1 \times S^1$, we consider the multivariate von Mises distribution, denoted by  $\mathcal{MVM}(\mu,\kappa, \Delta)$. The density at  $\theta \in \mathbb{T}$ is given by 

\begin{displaymath}
f(\theta; \mu,\kappa,\Delta)= \frac{1}{Z(\kappa,\Delta)} \exp \{  \kappa^\top c(\theta) +s(\theta)\Delta s(\theta) /2  \},
\end{displaymath}

where  $\mu \in \mathbb{T}^2$ (this parameter is called mean), $\kappa \geq 0 \in \mathbb{R}^d$  (concentration parameter),  $\Delta= (\lambda_{i,j})$ is a symmetric matrix on $\mathbb{R} ^{d \times d}$  with null diagonal entries ($\lambda_{i,i} = 0$ for all $i \in \{ 1,\ldots,d\}$), and  $Z(\kappa,\Delta)$ is a normalization constant. The functions $c_i$  and $s_i$ are defined by $c_i(\theta)= \cos(\theta_i - \mu_i)$ and $ s_i(\theta)= \sin(\theta_i - \mu_i)$ for all $i \in \{1, \ldots,d\}$.
In Figure \ref{toro} (left-hand panel), we show (in yellow) a sample of size 2000 from a $\mathcal{MVM}_1(\mu_1,\kappa_1, \Delta_1)$ distribution with

\begin{equation}\label{law1}
\mu_1=(\pi/2,0), \quad \kappa_1=(20,20), \quad \Delta_1= \left( \begin{matrix}  0&1 \\ 1 &0\end{matrix}\right). 
\end{equation}
In the right panel of Figure \ref{toro}, we show (in yellow) a sample of size 2000 from a mixture law given by

\begin{equation} \label{mixture}
0.4 \mathcal{MVM}_1(\mu_2,\kappa_1, \Delta_1) + 0.6 \mathcal{MVM}_2(\mu_3,\kappa_1, \Delta_1),
\end{equation}

where $\mu_2=(\pi/2,0)$ y  $\mu_3=(\pi/2,\pi/4)$ . 
In all cases, we consider $\lambda=0.8$ and bandwidth $h=0.2.$ The boundary of the theoretical level set is shown in red, while the boundary of the estimator is shown in magenta.

The Hausdorff distances between the theoretical curve and the estimated curve are $0.066$ and $0.107$.

\subsection{The half-sphere}

Finally, we considered the sphere $S^2 \subset \mathbb{R}^3$ endowed with the Riemannian metric inherited from $\mathbb{R}^3$. The sample is drawn from a the mixture of two von Mises--Fisher distributions given by 
\begin{equation}\label{von}
	f(x,\mu,\kappa)= C(\kappa) e^{\kappa \mu^\top x} \ind_{S^2}(x),
\end{equation}
 where $\kappa \geq 0$  and $\mu \in S^2$ are the concentration and directional mean parameters, respectively. $C(\kappa)$ is the normalizing constant; see \cite{mar72}.

\begin{figure}[h]
	\begin{center}
	\includegraphics[scale=.25]{ 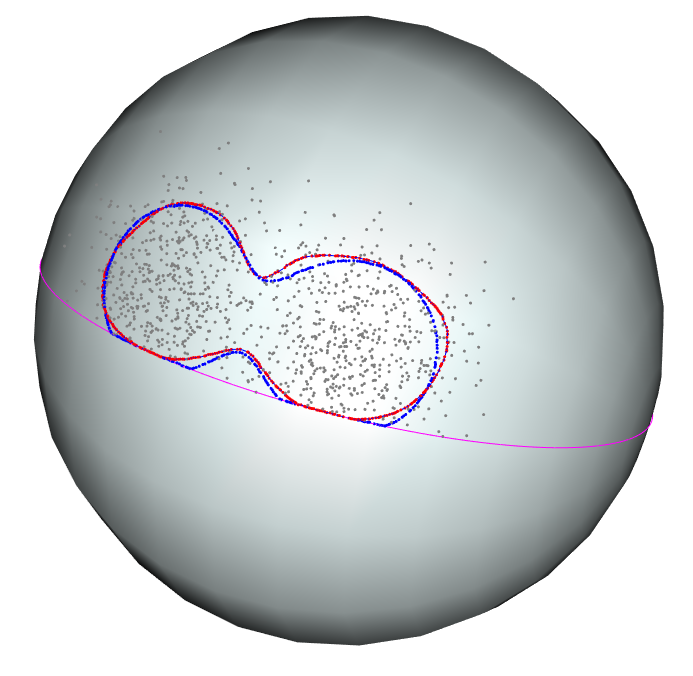}
		\includegraphics[scale=.4, height=5cm, width= 9 cm]{ 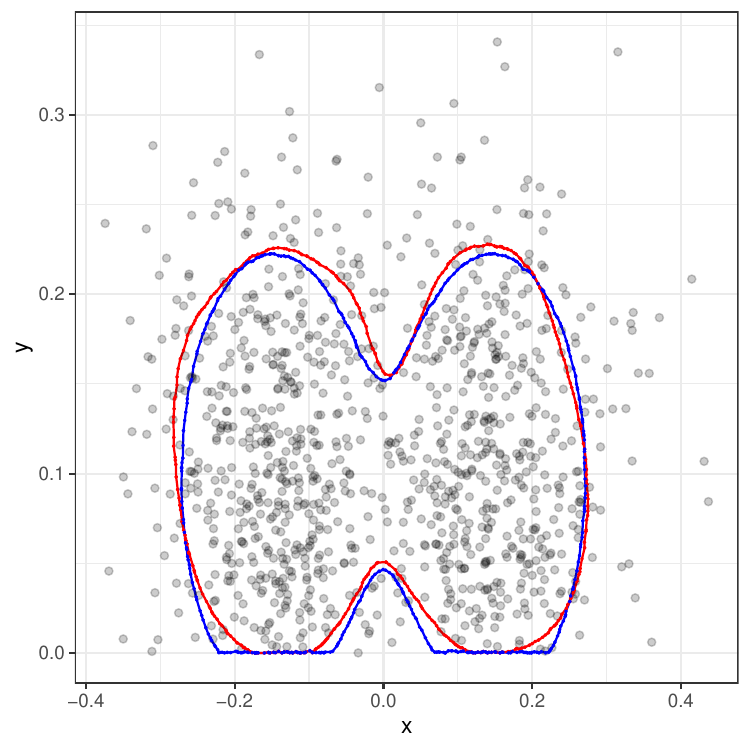}
		\caption{Left-hand panel: A sample of size $1000$ from the mixture of two von Mises--Fisher distributions given in equation \eqref{eqesf}. Right-hand panel: the stereographic projections of the sample and the level sets. In both cases, the estimator is shown in red, while the true underlying level set is shown in blue.}
		\label{esfera}
	\end{center}
\end{figure}

The mixture  is given by,
\begin{equation}\label{eqesf}
    f_0(\cdot)= 0.5 f \left( \cdot ,v_1/\Vert v_1 \Vert,40 \right) + 0.5 f\left(\cdot ,v_2/\Vert v_2 \Vert,40\right). 
\end{equation}
 with $v_1= (-1,-0.3,0.2)$ and $v_1= (-1,0.3,0.2)$.  Let $X=(X_1,X_2, X_3) \in S^2$  with distribution $f_0$, we consider the truncated random vector  $X^{\textrm{tr}}=X\mathbb{I}_{S^2_+}$, where $S^2_{+}=S^2 \cap \{ (x_1,x_2,x_3) \in \mathbb{R}^3 : x_3 \geq 0 \}$.
\\
In Figure \ref{esfera}, we show (left-hand panel) a sample of size $n=1000$ on $S^2_{+}$ of $X^{\textrm{tr}}$ , together with the estimated level set (in red ) and the true level set (in blue) at $\lambda=1.1$. In the right-hand panel, we show the stereographic projections of the sample and the estimators. The Hausdorff distance between the theoretical curve and the estimated (on the stereographic projections) curve is $0.019$.

\subsection{The two dimensional sphere}

To asses the performance of \eqref{nuestroest} and compare it with the density estimator proposed in \cite{hall} (which will be denoted by $\tilde{f}_n$),  specially designed for spherical data, we considered a Von-Mises distribution in $S^2$, whose density is given by \eqref{von}, with parameters  $\kappa=40$ and  $\mu=(0,0,1)$. Four sample sizes where considered: $n\in \{500,1000,5000,10000\}$ and the whole procedure was repeated 100 times.  For each replication the bandwidth $h$ for our estimator was selected as the one minimizing $\|\hat{f}_n-f\|_\infty:=\sup_{x\in S^2}|\hat{f}_n(x)-f(x)|$ on a grid of values, whereas the parameter $\kappa$ in \cite{hall}, was selected as the one minimizing $\|\tilde{f}_n-f\|_\infty$ on a grid of values. The mean over 100 replications is shown in Table \ref{comparacionhall}, while between $[]$ and $()$ we report the median and standard deviation respectively. The results are very encouraging. The performances of both estimators are very similar, with a slightly better behaviour of $\tilde f_n$ for $n=500, 1000$ while for $n=5000, 10000$ $\hat f_n$   outperforms $\tilde{f}_n$ slightly.

  \begin{table}[ht]
	\begin{center}
		\scalebox{1}{
			\begin{tabular}{ r| c|c}
				$n$   & $\|\hat{f}_n-f\|_\infty$ & $\|\tilde{f}_n-f\|_\infty$ \\ \hline
				500   &   1.1067  [1.0640] (0.2983)        &  1.1062  [1.0664] (0.2992)               \\
				1000  &   0.9651  [0.9596] (0.2079)        &  0.9631  [0.9601] (0.2092)              \\
				5000  &   0.6282  [0.6155] (0.1329)        &  0.6285  [0.6118] (0.1327)               \\
				10000 &   0.5262  [0.5018] (0.0910)        &  0.5358  [0.5164] (0.0915)
 			\end{tabular}}
		\caption{ Mean over 100 replications for  $\|\hat{f}_n-f\|_\infty$  and $\|\tilde{f}_n-f\|_\infty$. Between $[]$ and $()$ the median and standard deviation respectively.}
		\label{comparacionhall}
	\end{center}
\end{table}

\section{ A real data example: extreme and non-extreme winds in Uruguay}\label{winds}

As an example of real-data on a manifold, we will study the behavior of extreme and non-extreme winds, measured at the meteorological station located in the international airport of Carrasco, Montevideo Uruguay. The aim is to characterize the direction of these winds, with respect to the time period  at which the data were obtained  (divided in two groups, corresponding to cold or hot seasons respectively). The data are stored as points in a cylinder $\mathcal{C}$ (i.e, $\mathcal{C}=S^1\times \mathbb{R}^+$), with two parameters, the wind intensity $a\geq 0$ (in m/s), and the wind direction ($\theta\in [0,2\pi)$), see Figure \ref{cilindro}. The angles are measured clockwise, with zero located at North. The original database consisted of 149040 wind measures obtained between the dates 01/01/2000 and 31/12/2016. However, we ignored some missing data ($0.86\%$ in total) and we will study separately  extreme and   moderate (non-extreme winds). 
 
The intensities recorded are the average of the horizontal component of the wind in the last ten minutes of each hour, while the directions recorded are the average of the directions, also in the last ten minutes of each hour hour.  Gust are not recorded (sudden increases in the instantaneous wind,that  exceed the average wind in more than $ 5 $ m / s.)

 \subsection{Extreme winds}
 
\begin{figure}[p]
	\begin{center}
		\includegraphics[scale=.64]{ 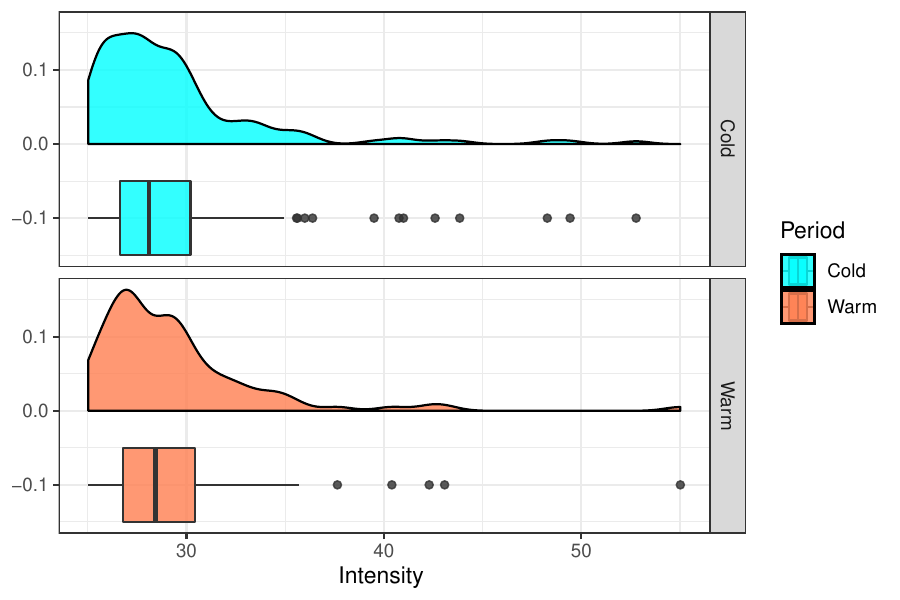}	
		\caption{Estimated marginal densities of extreme wind speed. Top panel: warm period. Bottom panel: cold period }
		\label{intensidades}
	\end{center}
\end{figure}

At large scale there exists two kinds of wind: the synoptic extreme wind, which are  produced by the passage of cold or warm fronts, with intensive convective activity, and the extreme non-synoptic winds, produced by extratropical cyclones associated with low and high pressure systems (see \cite{dura07}). Extreme winds are not classified according to the casuistry of the phenomenon.
According to \cite{dura12}, extreme maximum speeds are considered to be those that exceed $22$ m/s, because this is the threshold at which  different kind of risks for buildings and the population can exist. Measures below this threshold were removed. This is a typical procedure in the theory of extreme winds, and is called peak over threshold, (POT). We obtained 256 measure above this value, $39.8\%$ correspond to the warm period and $60.2\%$ to the cold period. The mean intensities are similar in both periods ($29.40$ and $29.38$ for the warm and cold period respectively).
In Figure \ref{intensidades}, we show the estimated marginal densities of the extreme wind speed in every period and their boxplots. As it can be observed, there is no evident significant difference between them.

\begin{figure}[p]
	\begin{center}
		\includegraphics[scale=.49]{ 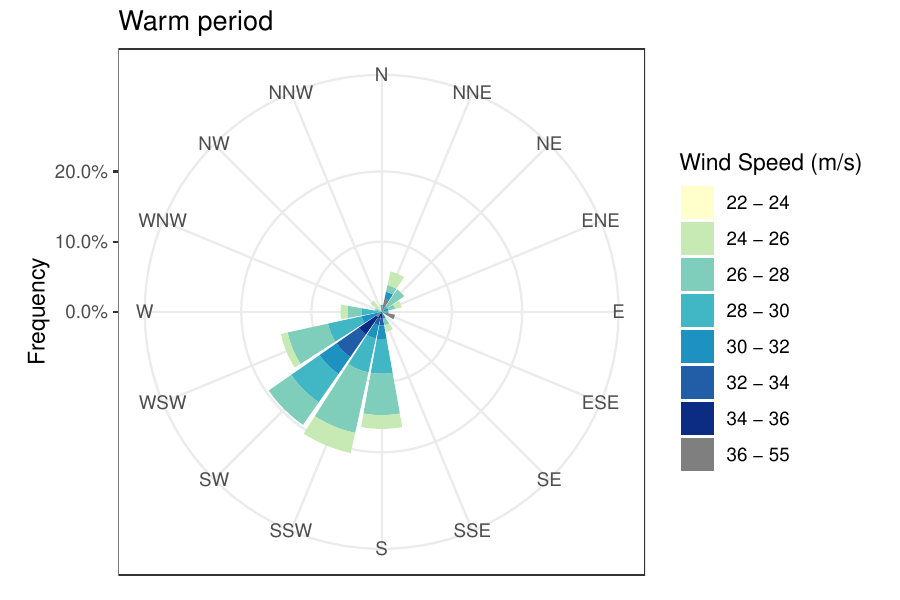}	
				\includegraphics[scale=.49]{ 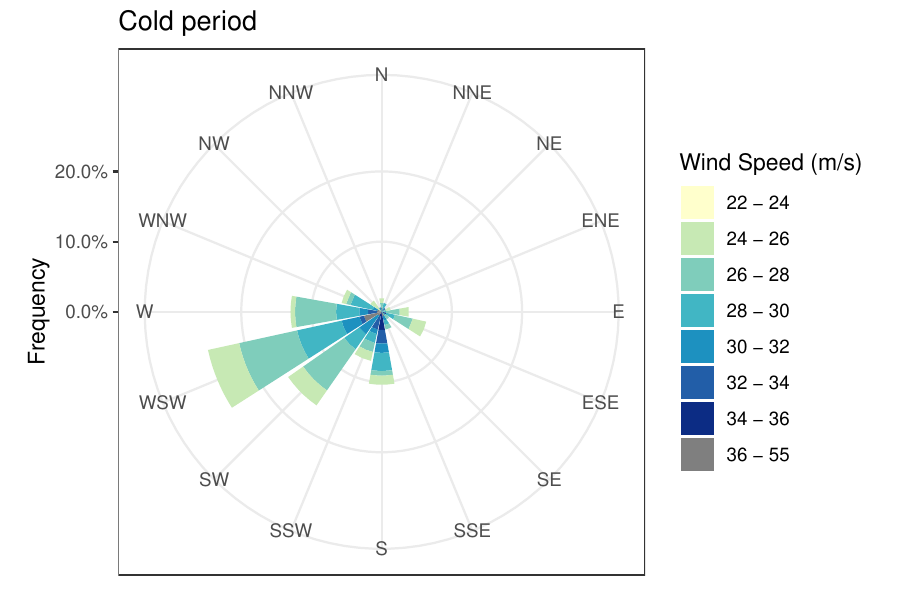}	
		\caption{Circular histograms. Left panel: warm period. Right panel: cold period. }
		\label{direcciones}
	\end{center}
\end{figure}

However, if we take into account the direction of the wind, the situation is different. Figure \ref{direcciones} shows the circular histograms of the marginal wind directions for the two periods of time considered. As it is seen, extreme winds in the warm period come more frequently from direction SW and SSW, with a strong trend to the south. At the cold period they come from the SW also, but with a trend to the West. This is also supported with the results obtained with the estimation of the level sets at each period.

\begin{figure}[p]
	\begin{center}
		\includegraphics[scale=.24]{ 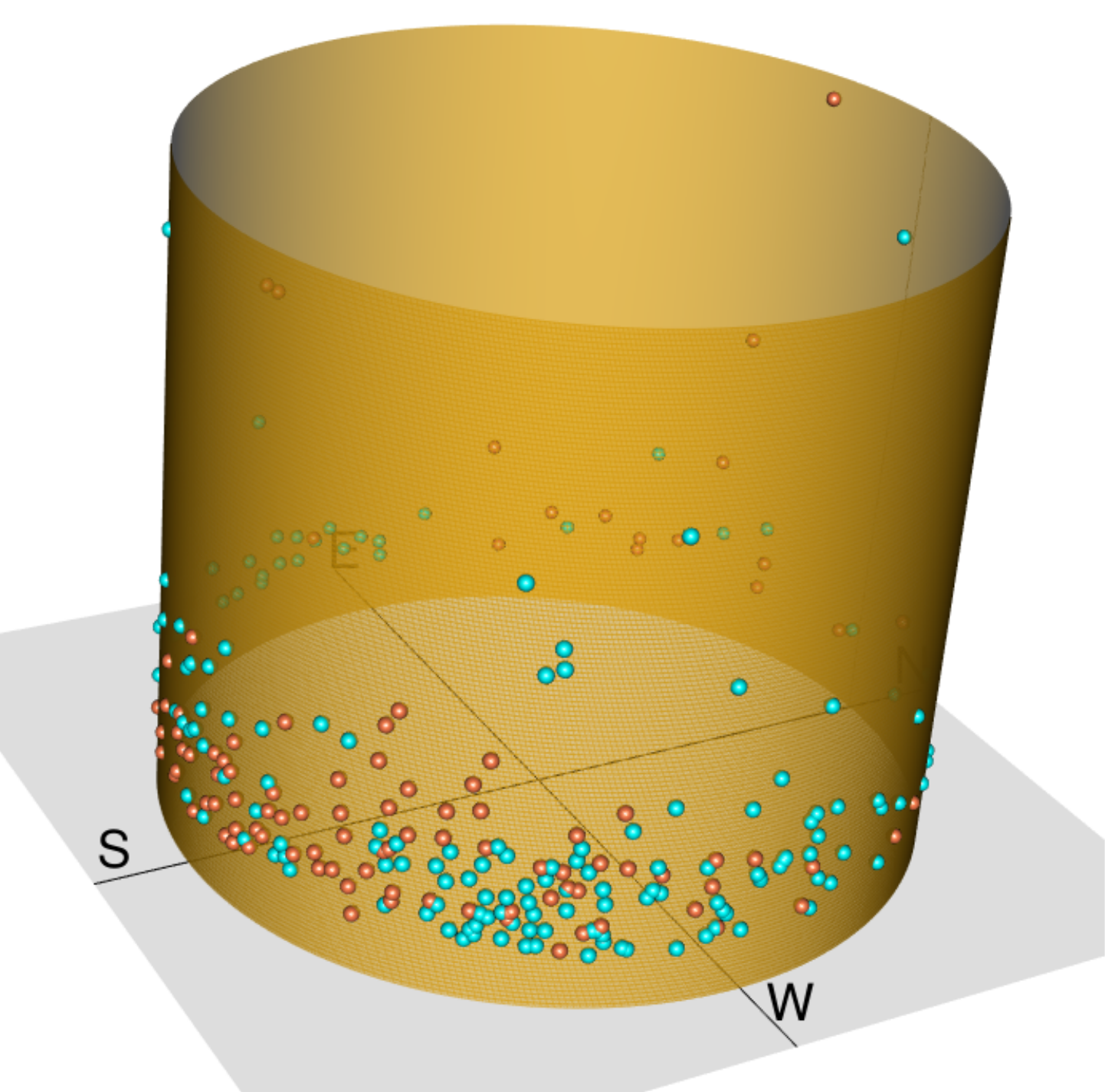}	
		\includegraphics[scale=.40]{ 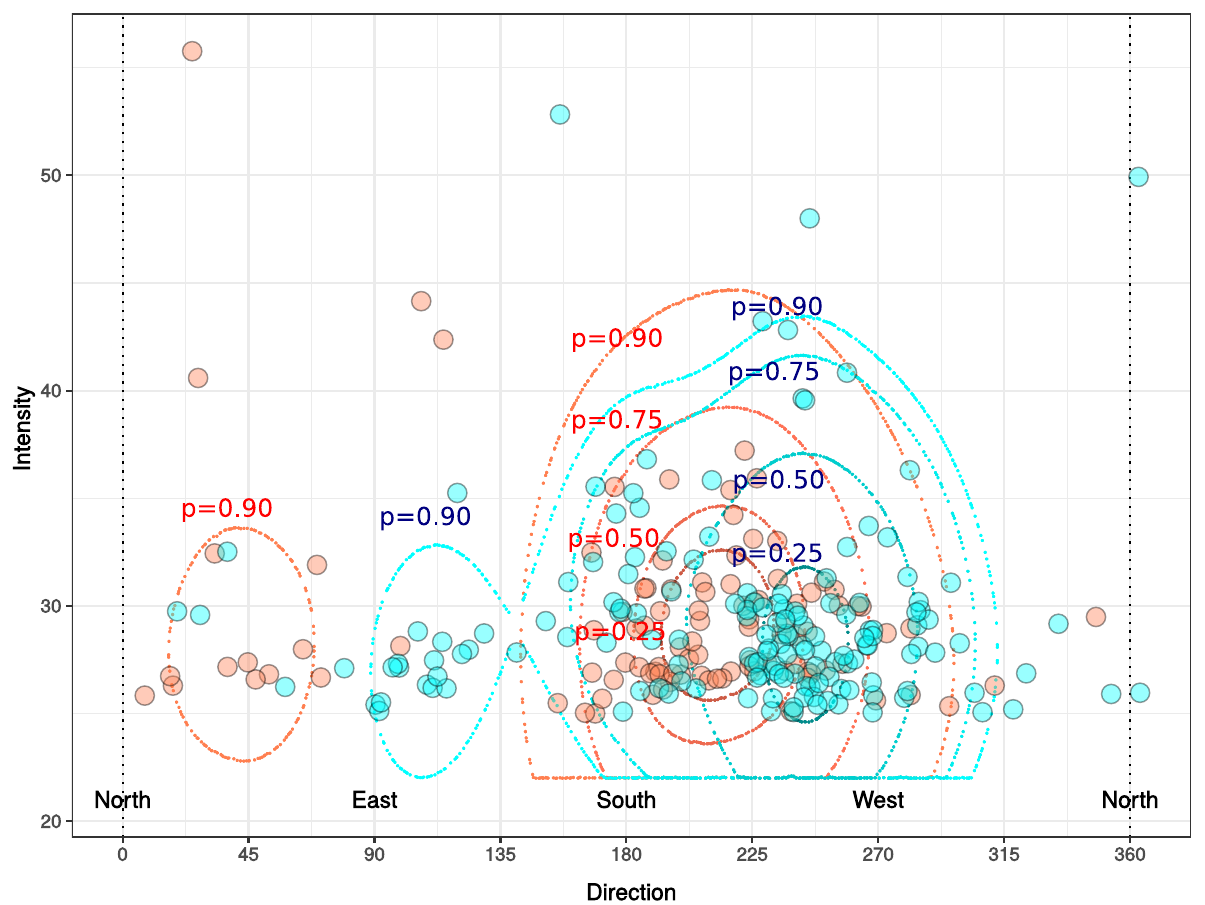}	
		\caption{Left-hand panel: the data on the cylinder $\mathcal{C}$, in sky-blue the data points corresponding to the cold period. Right-hand panel: the cylindrical coordinate map of data and the level sets for each period.  }
		\label{cilindro}
	\end{center}
\end{figure}

Figure  \ref{cilindro} shows at left the extreme winds on the cylinder and at right the estimation of four level sets for each period considered plot together, as it was proposed in Section \ref{densest}. 
The chosen levels where those containing $90 \%$, $75 \%$, $ 50 \%$ and $25 \%$ of the data.  As it can be seen most of the extreme winds in the cold period comes from SW and W, and most of the extreme winds in the warm period comes from S and SW.

As can be seen in Figure \ref{cilindro}, there is a shift in the direction of the extreme winds in the cold period with respect to the warm period in all level sets. In addition,  the shapes of the level sets are quite different. For the $90\%$ level set we found in both cases two connected components, which are more separated for the warm period than for the cold period. The small component in the warm period is centered around the North-East, while the small component in the cold period is centered around East.
In addition, the $10\%$ of the extreme winds with intensity above 40 m/s are located at very different directions in the cold period and in the warm period. In Table \ref{tab1} we provide the mean values for intensity and direction for each of the chosen level sets and for each period. 

\begin{figure}[p]
	\begin{center}
		\includegraphics[scale=.4]{ 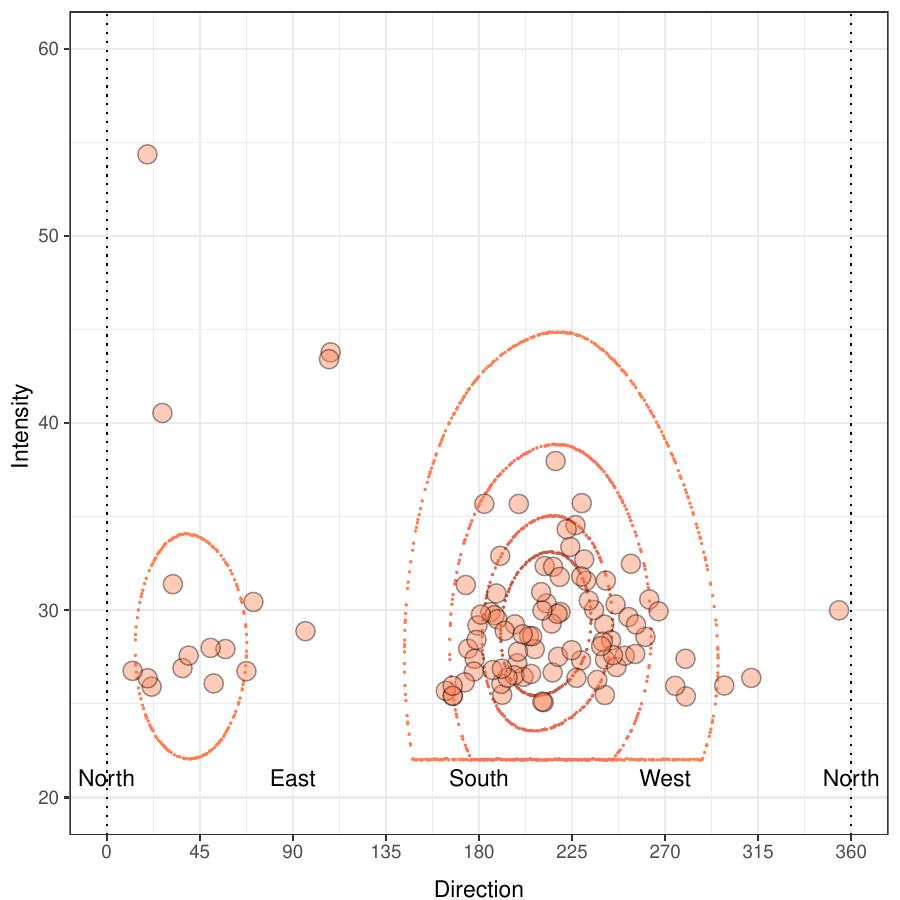}	
		\includegraphics[scale=.4]{ 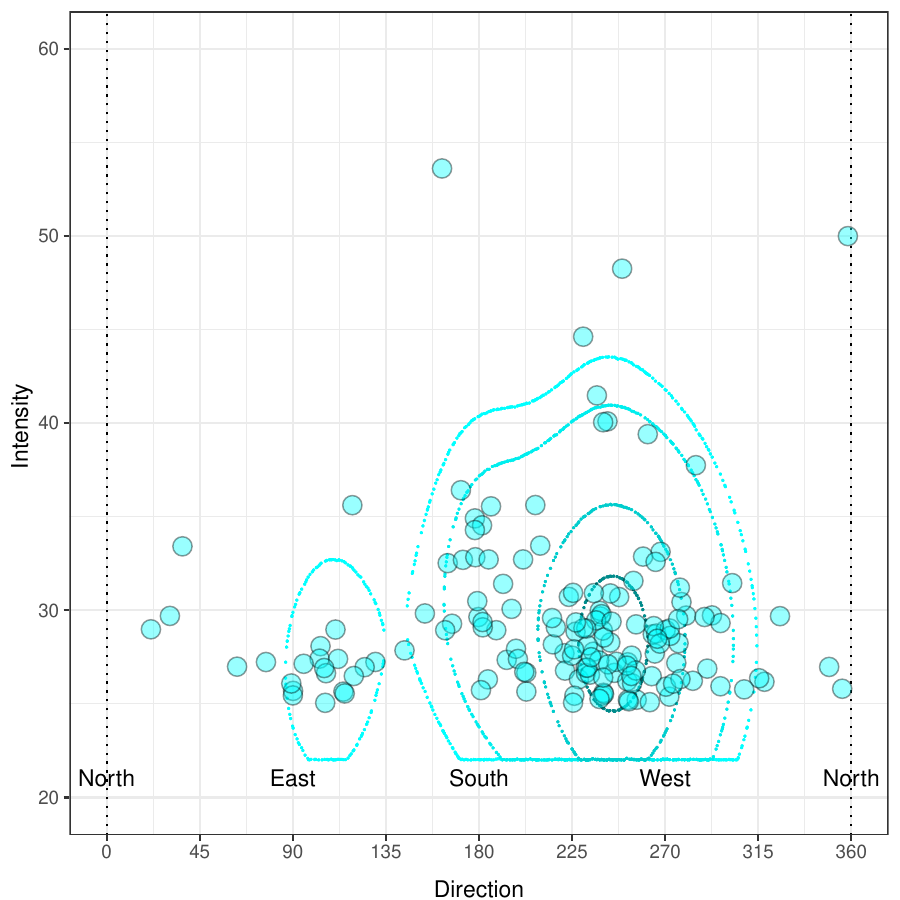}	
		\caption{Left-hand panel: the cylindrical coordinate map of data and the level sets for the warm period. Right-hand panel: the cylindrical coordinate map of data and the level sets for the cold period.}
		\label{cilindro2}
	\end{center}
\end{figure}

\begin{table}[p]\label{tab1}
	\begin{center}
		\begin{tabular}{c|c|c}
			& Direction & Intensity \\ 
			\hline
			$90\%$ &    200    &   28.92   \\
			&  (57.54)  &  (2.76)   \\
			$75\%$ &    214    &   29.29   \\
			&  (24.99)  &  (2.78)   \\
			$50\%$ &    214    &   29.00   \\
			&  (17.22)  &  (2.38)   \\
			$25\%$ &    212    &   29.04   \\
			&  (10.96)  &  (1.84)   \\ \hline
		\end{tabular}
		\hspace{2cm}
		\begin{tabular}{c|c|c}
			& Direction& Intensity  \\ 
			\hline
			$90\%$  & 220      &  28.88       \\
			& (53.14)  &  (3.28)\\
			$75\%$  &  236     &  28.93  \\
			&  (32.92) &  (3.06) \\
			$50\%$  &   247    &  27.99   \\
			&   (17.30)&  (1.94)\\
			$25\%$  &   241    &  27.65 \\
			& (8.23)   &  (1.66)\\
			\hline
		\end{tabular}
		
	\end{center}
	\caption{Mean value and standard deviation for the direction and the intensity, for each level set. Left: warm period. Right: Cold period.}
\end{table}

\begin{figure}[p]\label{f1}
	\begin{center}
		\includegraphics[scale=.35]{ 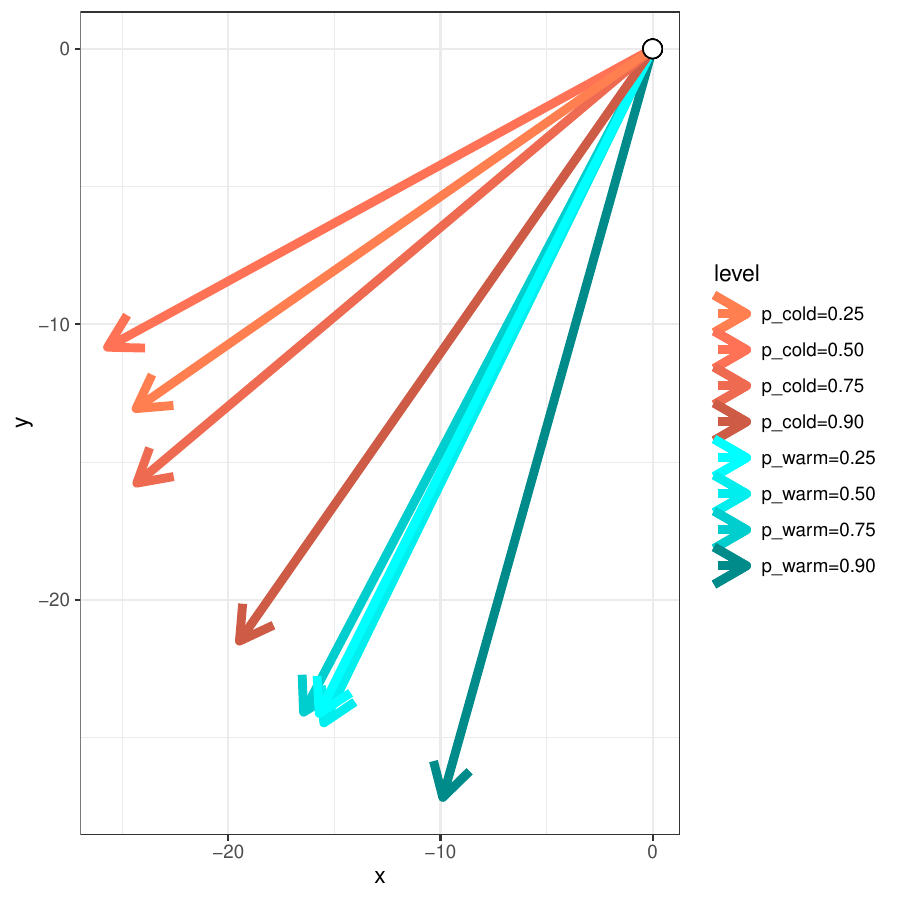}
		\caption{For each level we drawn an arrow whose norm is equal to mean intensity in that level, pointing toward the mean of the direction in that level. }	
	\end{center}
\end{figure}

\subsection{Moderate winds}
 Uruguay is currently one of the countries with the largest presence of wind energy in its electricity supply. Today almost 700 wind turbines are distributed in about thirty public and private parks. Uruguay occupies the second place in the world with 40.1\% of wind generation. Analyzing the behavior of moderate winds is important in many aspects in particular for wind energy generation.
In this subsection we replicate the previous analysis for the case of non-extreme winds, which corresponds to intensities below  22 m/s.

\begin{figure}[p]
	\begin{center}
		\includegraphics[scale=.450]{ 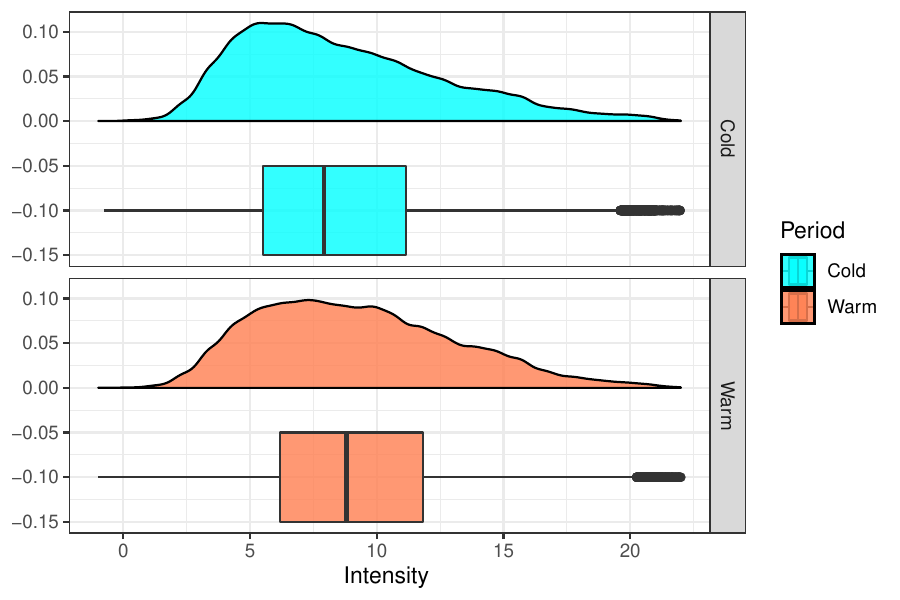}	
		\caption{Non-extreme winds. Estimated marginal densities of extreme wind speed. Top panel: warm period. Bottom panel: cold period.}
		\label{Ncil1}
	\end{center}
\end{figure}

\begin{figure}[p]
	\begin{center}
		\includegraphics[scale=.450]{ 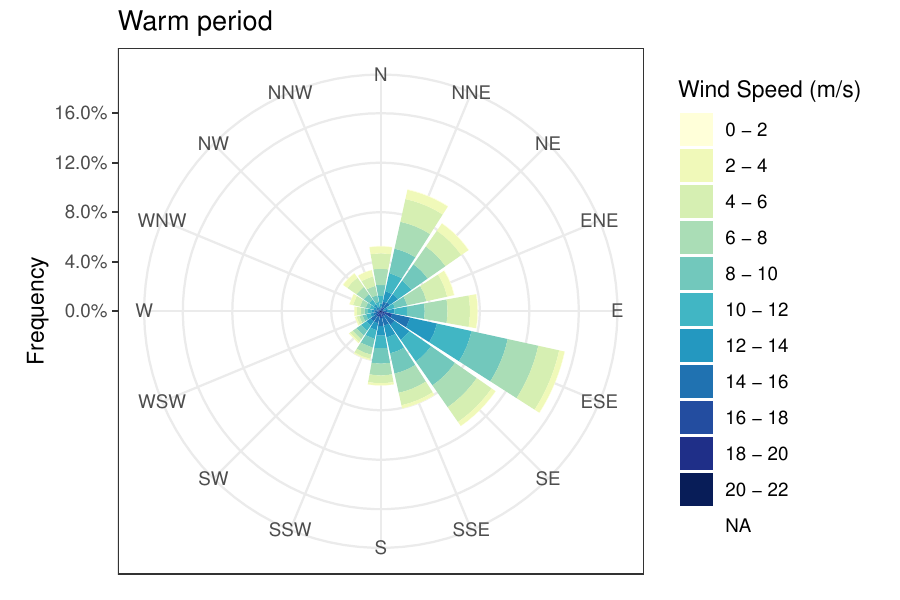}	
		\includegraphics[scale=.450]{ 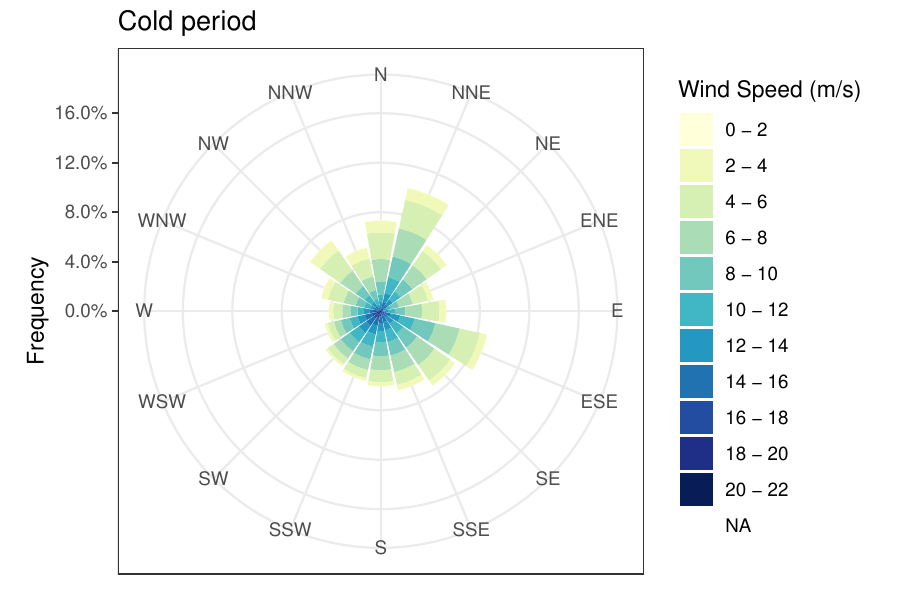}	
		\caption{Circular histograms for non-extreme winds. Left-hand panel: warm period. Right-hand panel: cold period. }
		\label{Ncil2}
	\end{center}
\end{figure}

In Figure \ref{Ncil1} we show the estimated marginal densities of the non-extreme wind speed at every period and their boxplots, as  can be seen they are quite similar.  Figure \ref{Ncil2} shows the circular histograms of the marginal wind directions for the two periods of time considered. As can be seen, the directions of non-extreme winds in the warm period corresponds more frequently with ESE direction, with a strong trend to the south. During cold period the main directions corresponds to  NNE. This is also supported by the results obtained with the estimation of the level sets at each period.

\begin{figure}[p]
	\begin{center}
		\includegraphics[scale=.450]{ 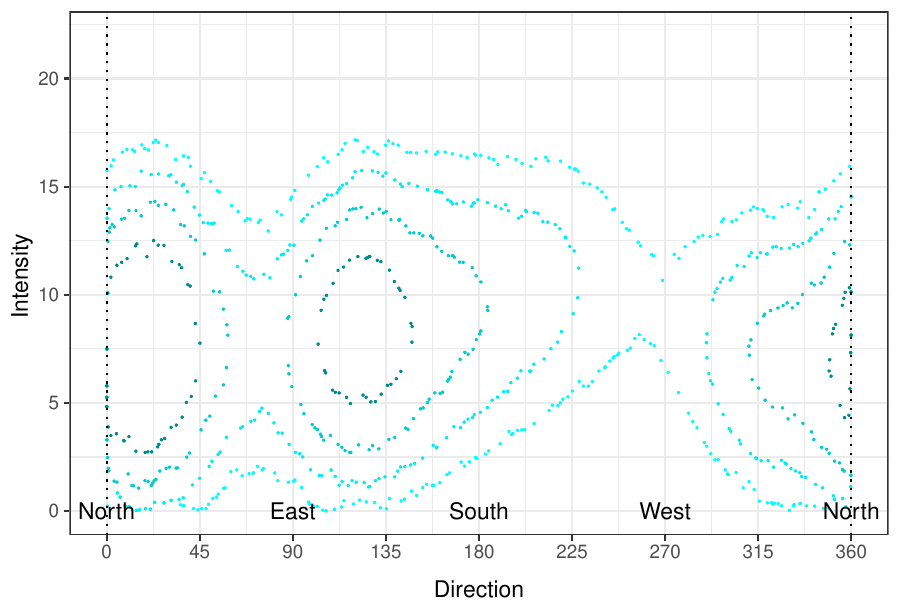}	
		\includegraphics[scale=.450]{ 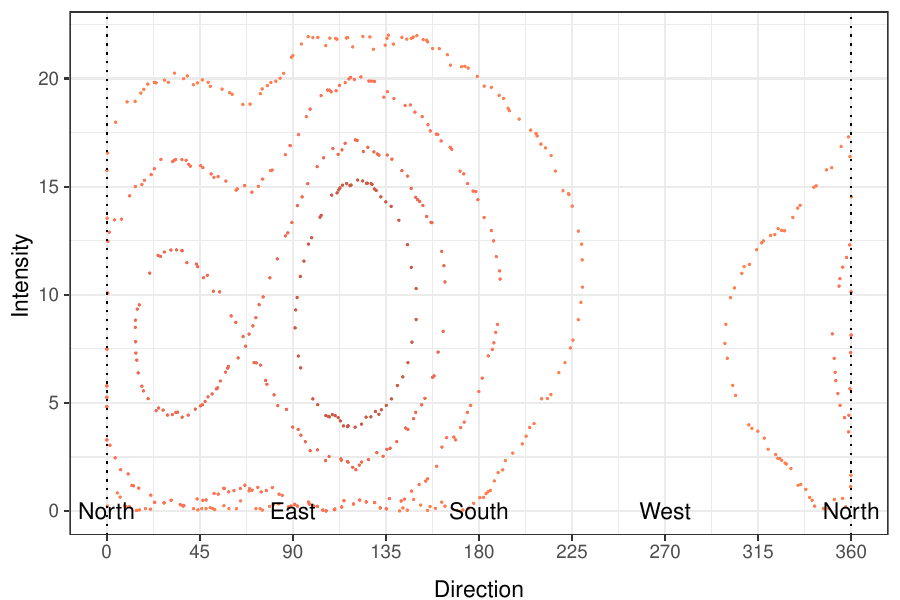}	
		\caption{Non-extreme winds. Left-hand panel: the level four sets for cold period. Right-hand panel: the four level sets for the warm period.}
		\label{Ncil3}
	\end{center}
\end{figure}

\begin{figure}[p]\label{f2}
	\begin{center}
		\includegraphics[scale=.35]{ 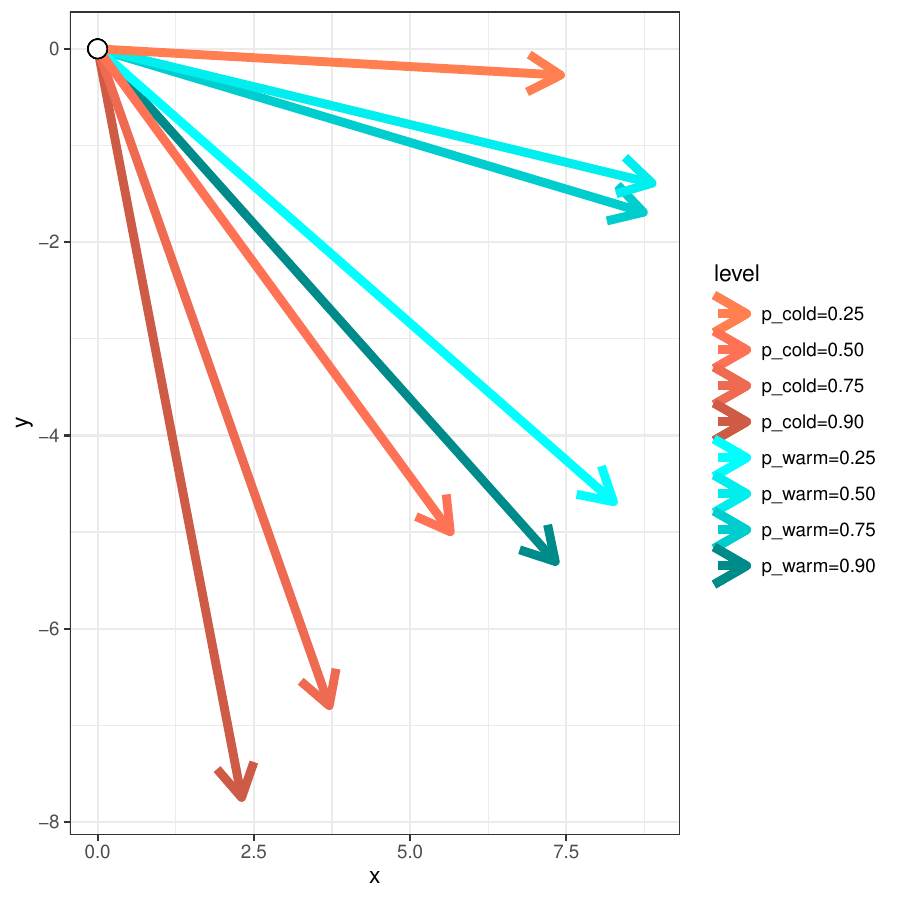}
		\caption{Non-extreme winds.For each level we drawn an arrow whose norm is equal to mean intensity in that level, pointing toward the mean of the direction in that level. }	
	\end{center}
\end{figure}

 Figure  \ref{Ncil3} shows  the four level sets for each period considered, as proposed in Section \ref{densest}. 
The chosen levels where those containing $90 \%$, $75 \%$, $ 50 \%$ and $25 \%$ of the data.  As can be seen the directions of most of the non-extreme winds in the cold period corresponds to SW and W, while those of the extreme winds in the warm period corresponds to directions S and SW (see Figure \ref{f2} and Table \ref{tab2}).

\begin{table}[t]\label{tab2}
	\small
	\begin{center}
		\begin{tabular}{c|c|c}
			       & Direction & Intensity \\ \hline
			$90\%$ &    126    &   9.04   \\
			       &  (87.92)  &  (3.72)   \\
			$75\%$ &    100    &   8.89   \\
			       &  (63.85)  &  (3.50)   \\
			$50\%$ &    99     &   8.98   \\
			       &  (40.83)  &  (2.96)   \\
			$25\%$ &    119    &   9.49   \\
			       &  (13.94)  &  (2.64)   \\ \hline
		\end{tabular}
		\hspace{2cm}
		\begin{tabular}{c|c|c}
			       & Direction & Intensity \\ \hline
			$90\%$ &    163    &   8.08   \\
			       &  (109.01)    &  (3.27)   \\
			$75\%$ &    151    &   (7.74)  \\
			       &  (111.44)  &  (2.95)   \\
			$50\%$ &    131    &   7.53   \\
			       &  (113.19)  &  (2.59)   \\
			$25\%$ &    92    &   7.41   \\
			       &  (102.48)   &  (2.12)   \\ \hline
		\end{tabular}
		
	\end{center}
	\caption{Non-extreme winds. Mean value and standard deviation for the direction and the intensity, for each level set. Left: warm period. Right: Cold period.}
\end{table}

\newpage

\section*{Appendix}

\textit{Proof of Proposition \ref{unitang}}
	
	From Proposition 14 in \cite{tha08} $\partial M$ and $M$ have positive reach. Denote by $\text{reach}(M)$ the reach of $M$, (it is also known as the condition number, see \cite{gen12b}). By point 4 in Lemma 3 , in \cite{gen12b} $R(x,y)=\|x-y\|/\rho(x,y)$ is bounded away from zero if $\|x-y\|<\text{reach}(M)/2$. Since $\text{reach}(\partial M)>0$ w.r.t. to the intrinsic metric in $M$, (see remark 12 in \cite{tha08}),  $M$ admits a normal collard $N_r$. Since $M\setminus \exp(N_r)$ is compact and of class $\C^2$ the injectivity radius on $M\setminus \exp(N_r)$ is bounded away from zero, from where it follows that $M$ is uniformly tangible.

\textit{Proof of Theorem \ref{convunif2}}

Let us bound
\begin{multline*}
\sup_{x\in M_n}|\hat{f}_{n,h}(x)-f(x)|\leq \sup_{x\in M_n} |\hat{f}_{n,h}(x)-f_{n,h}(x)|+\\
\sup_{x\in M_n}|f_{n,h}(x)-E(f_{n,h}(x))|+
\sup_{x\in M_n}|E(f_{n,h}(x))-f(x)|=A_1+A_2+A_3
\end{multline*}

First we prove that $h^{-(1+\delta)}A_2\to 0 $ a.s., for all $\delta<1/2$.  We will follow the same ideas used in \cite{hr09}. First, we define the random variables 
$$V_j(x)=\frac{1}{m_0(x)} K\Big(\frac{\|x-X_j\|}{h}\Big)- \frac{1}{m_0(x)}E\Big[K\Big(\frac{\|x-X_j\|}{h}\Big)\Big]$$
and let $S_n(x)=\sum_{j=1}^n V_j(x)$. Observe that from \eqref{boundm0}, $m_0(x)>1/2$ for all $x$ and $n$ large enough (independent of $x$). Because $K$ is bounded, it follows that $|V_j(x)|\leq C_2$ for all $x$. Let $\beta_n=h^{-(1+\delta)}\to \infty$, where $\delta<1/2$, then from Bernstein's inequality,
\begin{equation}\label{b2}
\sup_{x\in M_n} P\Big(\beta_n\frac{1}{nh^{d'}} |S_n(x)|>\epsilon\Big)\leq 2\exp\Big(-\frac{ \epsilon^2}{4C_2^2}\frac{n h^{d'}}{\beta_n^2}\Big)
\end{equation}
For $n$ large enough fixed, we consider a finite collection of balls $B_i=B(p_i,h^{\gamma})$ centered  at $p_i\in M_n$, with  $\gamma>d'+1$ such that $\beta_n h^{\gamma-d'-1}\to 0$, and $M_n\subset \cup_{i=1}^l B_i$. Because $M_n$ is compact and $\gamma>d'+1$, $l\leq C_3 h^{-\gamma}$.
$$\sup_{x\in M_n} \frac{1}{nh^{d'}}|S_n(x)|\leq \max_{1\leq j\leq l} \sup_{p\in B_j} \frac{1}{nh^{d'}}|S_n(x)-S_n(p_j)|+\max_{1\leq j\leq l}\frac{1}{nh^{d'}}|S_n(p_j)|=I_1+I_2$$
Because $K$ is Lipschitz, $\beta_n I_1\leq C_4 \beta_n h^{\gamma-(d'+1)}\to 0$ for some positive constant $C_4$. Then $I_1<\epsilon/2$ for $n$ large enough. From \eqref{b2} we get that for $n$ large enough
$$P\Big(\beta_n\sup_{x\in M_n} \frac{1}{nh^{d'}}|S_n(x)|>\epsilon\Big)\leq P\Big(\beta_n I_2 >\frac{\epsilon}{2}\Big)\leq 2C_3h^{-\gamma}\exp\Big(-\frac{C_5n h^{d'}}{\beta_n^2}\Big),$$
$C_5$ being a positive constant. Now from Borel-Cantelli's lemma, together with condition $nh^{d'}/(\beta_n^2\log(n))\to \infty$, it follows that $\beta_n A_2\to 0$ $a.s$.

To bound $A_3$, first we use that the term $\mathcal{O}_x(h^2)$ can be bounded independently of $x$, from above by $C_1h^2$ for some constant $C_1>0$ (see the proof of Theorem 3.1 in \cite{bs17}). Let $c_n=\inf_{x\in M_n} \rho(x,\partial M)$.  Now if we bound $m_1(x)\leq (1/(2\sqrt{\pi}))\exp(-c_n^2/h^2)$ and using that $\|\nabla f(x)\|\leq C$, then it follows from \eqref{bias} that $A_3$ is of the order $o(h)$.

To bound $A_1$ observe that   $\hat{f}_{h,n}(x)=m_0(x)f_{n,h}(x)$,

$$\sup_{x\in M_n} |\hat{f}_{n,h}(x)-f_{n,h}(x)|\leq \sup_{x\in M_n}|m_0(x)-1|\sup_{x\in M_n} f_{n,h}(x)$$

From $b_x\geq c_n>0$ for all $x\in M_n$ it follows that 
\begin{equation}\label{boundm0}
|m_0(x)-1|= \frac{1}{\sqrt{\pi}}\int_{\frac{b_x}{h}}^{+\infty} \exp(-z^2)dz\leq \frac{1}{\sqrt{\pi}}\int_{\frac{c_n}{h}}^{+\infty} \exp(-z^2)dz\leq \frac{h}{2c_n\sqrt{\pi}}\exp\Big(-\frac{c_n^2}{h^2}\Big)
\end{equation}

To bound $\sup_{x\in M_n} f_{n,h}(x)$, we proceed as we did with $A_3$, and it follows that 
$$\beta_n\sup_{x\in M_n} |f_{n,h}(x)-E(f_{n,h}(x))|\to 0\quad a.s.$$
so it is enough to bound $|E(f_{n,h}(x))|$, but we have proven that $\sup_{x\in M_n} |E(f_{n,h}(x))-f(x)|\to 0$. Then, because $f$ is continuous and $M$ is compact, it is bounded. So, for $n$ to be large enough, $\sup_{x\in M_n} |E(f_{n,h}(x))|<2\sup_{x\in M} f(x)<\infty$. Finally, we have proven that 
$$A_1=\sup_{x\in M_n} |\hat{f}_{n,h}(x)-f_{n,h}(x)|=o(h/c_n).$$\\

\textit{Proof of Theorem \ref{th:3}}\\
Let us prove point 1, because $  L_f(\lambda)\cap \partial M=\emptyset$, we can take $\delta>0$ small enough such that $L_f(\lambda)\subset M_{\delta}$, then condition $M1$ in \cite{cgmrc:06} is fulfilled in $M_\delta$; that is,  $B_\rho(x,r)$ is connected for all $x\in M_\delta$ and for all $0<r<\delta$. Because $M$ is compact, condition f2  in \cite{cgmrc:06} is fulfilled. Then,  Theorem  \ref{convunif}, and Theorem 2 in \cite{cgmrc:06} entails that
$d_H(\partial L_{\hat{f}_{n,h}}(\lambda),\partial L_{f}(\lambda))\to 0$.
To prove 2, observe that Theorem 2.1 in \cite{mol98} implies that  $d_H(L_{\hat{f}_{n,h}}(\lambda),L_f(\lambda))\to 0$ (observe that $\partial L_f(\lambda)=\{x:f(x)=\lambda\}$). Finally, to prove point 3, observe that if $\nabla_x f\neq 0$ for all $x:f(x)=\lambda$,  then $\partial L_f(\lambda)=\{x:f(x)=\lambda\}$   is a $d'$-1,dimensional submanifold of $M$, and then $\mu(\partial L_f(\lambda))=0$, then point 3 is a consequence of Theorem 2 in \cite{cue12}, which still holds for any metric space.\\

\textit{Proof of Theorem \ref{convhaus}}

Let $\epsilon_n\to 0$ such that $\epsilon_n/h\to \infty$ and define the sequence of sets $M_{\epsilon_n}=\{x\in M:\rho(x,\partial M)\geq \epsilon_n\}$, observe that $M_{\epsilon_n}$ is compact for all $n$, then
\begin{multline*}\label{ineq}
d_H\big(L_{\hat{f}_{n,h}}(\lambda),L_{f}(\lambda)\big)\leq d_H\big(L_{\hat{f}_{n,h}}(\lambda),L_{\hat{f}_{n,h}}(\lambda)\cap M_{\epsilon_n}\big)+\\
d_H\big(L_{\hat{f}_{n,h}}(\lambda)\cap M_{\epsilon_n},L_{f}(\lambda)\cap M_{\epsilon_n}\big)+d_H\big(L_{f}(\lambda)\cap M_{\epsilon_n},L_{f}(\lambda)\big)=I_1+I_2+I_3.
\end{multline*}
To prove that $I_2\to 0$ a.s. let us denote $\gamma_n=\sup_{x\in M_{\epsilon_n}}|\hat{f}_{h,n}(x)-f(x)|$, then 
$$L_f(\lambda+\gamma_n)\cap M_{\epsilon_n}\subset L_{\hat{f}_{h,n}}(\lambda)\cap M_{\epsilon_n}\subset L_f(\lambda-\gamma_n)\cap M_{\epsilon_n}.$$
Then $I_2\leq \sup_{x\in L_f(\lambda-\gamma_n)\cap M_{\epsilon_n}}\rho(x,L_f(\lambda+\gamma_n)\cap M_{\epsilon_n})=:R_n$. To prove that $R_n\to 0$ assume by contradiction that there exists $\delta>0$ and $x_n\in L_f(\lambda-\gamma_n)\cap M_{\epsilon_n}$ such that $\delta<\rho(x,L_f(\lambda+\gamma_n)\cap M_{\epsilon_n})$. We can assume that $x_n\to x_0$ and  for all $n$, $\delta/2\leq \rho(x_0,L_f(\lambda+\gamma_n)\cap M_{\epsilon_n})$. If $f(x_0)>\lambda$ there exists $N_{x_0}$ such that $f(z)>\lambda$ for all $z\in N_{x_0}$,  fix $\epsilon<\delta/2$ and $z\in N_{x_0}$ with $\rho(x_0,z)<\epsilon$ and $n$ large enough such that $z\in M_{\epsilon_n}$, then for $n$ large enough  $z\in L_f(\lambda+\gamma_n)$, which is a contradiction. Then, $f(x_0)=\lambda$. Fix $a_j$ with $\rho(x_0,a_j)<\delta/2$ and $f(a_j)>\lambda$, then for $n$ large enough $f(a_j)>\lambda+\gamma_n$ and $a_j\in M_{\epsilon_n}$ which is again a contradiction. This proves that $R_n\to 0$ and then $I_2\to 0$.\\

Let us prove that $I_1\to 0$ a.s., as $n\to \infty$. Because $L_{\hat{f}_{n,h}}(\lambda)\cap M_{\epsilon_n}\subset L_{\hat{f}_{n,h}}(\lambda)$, it follows that,

$$d_H\big(L_{\hat{f}_{n,h}}(\lambda),L_{\hat{f}_{n,h}}(\lambda)\cap M_{\epsilon_n}\big)=\sup_{x\in L_{\hat{f}_{n,h}}(\lambda)} \rho(x,L_{\hat{f}_{n,h}}(\lambda)\cap M_{\epsilon_n}).$$

Suppose by contradiction that $I_1$ does not converge to $0$ a.s., then there exists $\delta>0$ such that with positive probability there exists $x_n\in L_{\hat{f}_{n,h}}(\lambda)$ such that $\rho(x_n,L_{\hat{f}_{n,h}}(\lambda)\cap M_{\epsilon_n})>\delta$. Because $M$ is compact, we can assume that there exists $x\in  M$ such that $x_n\to x$ (by taking a subsequence if it is necessary), with positive probability. Observe that $x\in \partial M$ because $x_n \in M\setminus M_{\epsilon_n}$ and $\epsilon_n\to 0$. Now let us prove that $f(x)\leq \lambda$ (recall that $\hat{f}_{h,n}(x_n)\geq \lambda$). Suppose that $f(x)> \lambda$,
then there exists a neighbourhood $B_x$ on $M$, of $x$, such that $f(z)>\lambda$ for all $z\in B_x$. Observe that for all $n$ large enough $B_x\cap M_{\epsilon_n}\neq \emptyset$, then for all $y\in B_x\cap M_{\epsilon_n}$, $\hat{f}_{h,n}(y)\geq \lambda$ because $\sup_{x\in M_{\epsilon_n}}\|\hat{f}_{h,n}(x)-f(x)\|\to 0$ a.s. But then we can choose a sequence $y_n\in  B_x\cap M_{\epsilon_n}$ with $y_n\to x$ such that $\hat{f}_{h,n}(y_n)\geq \lambda$. So $\rho(x_n,L_{\hat{f}_{n,h}}(\lambda)\cap M_{\epsilon_n})\leq \rho(x_n,y_n)\leq \rho(x_n,x)+\rho(x,y_n)\to 0$, which contradict that $\delta<\rho(x_n,L_{\hat{f}_{n,h}}(\lambda)\cap M_{\epsilon_n})$. This proves that $f(x)\leq \lambda$. 


Observe that  $\hat{f}_{h,n}(x_n)=m_0(x_n)f_{h,n}(x_n)\geq \lambda$. We will prove that $f_{h,n}(x_n)\to f(x)\leq \lambda$ a.s., which is a contradiction because $m_0(x_n)\to 1/2$. Let us bound,
$$|f_{h,n}(x_n)-f(x)|\leq |f_{h,n}(x_n)-E(f_{h,n}(x_n))|+|E(f_{h,n}(x_n))-f(x_n)|+|f(x_n)-f(x)|$$
and $|f_{h,n}(x_n)-E(f_{h,n}(x_n))|\leq \sup_{x\in M}|f_{h,n}(x)-E(f_{h,n}(x))|$. Now the convergence  $\sup_{x\in M}|f_{h,n}(x)-E(f_{h,n}(x))|\to 0$ a.s., is proved following the same ideas used to prove that $A_2\to 0$ as in the proof of Theorem \ref{convunif}). Regarding the second term, we also bound 
$|E(f_{h,n}(x_n))-f(x_n)|\leq \sup_{x\in M} |E(f_{h,n}(x))-f(x)|$, which converges to 0 following the same ideas used to prove $A_3\to 0$ in Theorem \ref{convunif}.
Finally, $f(x_n)\to f(x)$ because $f$ is a continuous function.\\

To prove that $I_3\to 0$ as $n\to \infty$ assume by contradiction that this is not true, then there exists $\delta>0$ and a sequence $x_n$ such that $f(x_n)\geq \lambda$, $x_n \in M\setminus M_{\epsilon_n}$ and $\rho(x_n, M_{\epsilon_n}\cap L_f(\lambda))>\delta$. Because $\epsilon_n\to 0$, there exists a subsequence of $x_n$ (which will be denoted $x_n$ for ease of writing), such that $x_n\to x\in \partial M$. Because $f$ is continuous $f(x)\geq \lambda$. If $f(x)>\lambda$, then there exists $N_x$ a neighborhood of $x$ such that  for all $y\in N_x\cap M_{\epsilon_n}$ $f(y)\geq \lambda$. Let us choose $y_n\to x$ and $y_n \in N_x\cap M_{\epsilon_n}$, then $\delta<\rho(x_n, M_{\epsilon_n}\cap L_f(\lambda))\leq \rho(x_n,y_n)\to 0$, which is a contradiction.  The other case is $f(x)=\lambda$, let $a_j\to x$ such that $f(a_j)>\lambda$ for all $j$. For all $j$, we can choose $n(j)\to \infty$ as $j\to \infty$, such that $a_j\in M_{\epsilon_{n(j)}}$. Then $$\delta<\rho(x_{n(j)}, M_{\epsilon_{n(j)}}\cap L_f(\lambda))\leq \rho(x_{n(j)},a_{n(j)})\leq \rho(x_{n(j)},x)+\rho(x,a_{n(j)})\to 0\text{ as } j\to \infty.$$  

\textit{Proof of Theorem \ref{th:5}}

Let us denote $\mathcal{X}_n=\{X_1,\dots,X_n\}$, then  
\begin{multline*}
d_H\big(C_r(\{X_i:\hat{f}_{h,n}(X_i)>\lambda\}),C_r(\{X_i:f(X_i)>\lambda\})\big)\leq \\
 d_H\big(C_r(\mathcal{X}_n\cap L_{\hat{f}_{h,n}}(\lambda)),L_{\hat{f}_{h,n}}(\lambda)    \big)+\\ d_H\big(L_{\hat{f}_{h,n}}(\lambda), L_f(\lambda)\big)+ d_H\big( L_f(\lambda),C_r(\mathcal{X}_n\cap L_{{f}}(\lambda))\big)=A+B+C.
\end{multline*}	
From Theorem \ref{convhaus}, $B\to 0$ a.s. Because $L_f(\lambda)$ is $r$-convex, $\mathcal{X}_n\cap L_f(\lambda)\subset C_r(\mathcal{X}_n\cap L_f(\lambda))\subset L_f(\lambda)$ and then,
$$C\leq d_H(\mathcal{X}_n\cap L_f(\lambda), L_f(\lambda))\to 0 \quad a.s., \text{ as } n\to \infty.$$
Regarding $A$,  observe that  $A=\sup_{x\in L_{\hat{f}_{h,n}(\lambda)}} \rho(x,C_r(\mathcal{X}_n\cap L_{\hat{f}_{h,n}}(\lambda)))$.  Let us proceed by contradiction, assume that with positive probability $A$ does not converge to $0$, then there exists a sequence $x_n \in  L_{\hat{f}_{h,n}}(\lambda)$ and $\delta>0$ such that $\delta<\rho(x_n,C_r(\mathcal{X}_n\cap L_{f}(\lambda)))$ for all $n>n_0$.  Because $M$ is compact, there exists a convergent subsequence of $x_n$ (which we will denote $x_n$) such that $x_n\to x$. Because $B\to 0$, it follows that $f(x)\geq \lambda$ but  with positive probability $\delta/2<\rho(x,C_r(\mathcal{X}_n\cap L_{\hat{f}_{h,n}}(\lambda)))$ for all $n$ large enough. If $f(x)>\lambda$, then there exists $\eta>0$ such that for all $z\in B_\rho(x,\eta)$ $f(z)>\lambda$. Let us take $0<\eta<\delta/2$, then with probability one, for $n$ large enough $\hat{f}_{h,n}(z)>\lambda$ for all $z\ B(x,\eta)$. Let us take $n$ large enough such that $d_H(\mathcal{X}_n\cap L_f(\lambda),L_f(\lambda))<\eta$, then $\mathcal{X}_n\cap B(x,\eta)\neq \emptyset$ but then $\rho(x,C_r(\mathcal{X}_n\cap L_{\hat{f}_{h,n}}(\lambda)))<\eta$ which is a contradiction. The case $f(x)=\lambda$ is proved in the same way, let $a_j$ such that $f(a_j)>\lambda$ and $0<\eta<\delta$ such that for all $z\in B_\rho(a_j,\eta)$ $f(z)>\lambda$. Let $n$ be large enough such that $\hat{f}_{h,n}(z)>\eta$ for all $z\in B(a,\eta)$ and $d_H(\mathcal{X}_n\cap L_f(\lambda),L_f(\lambda))<\eta$.  Again $\rho(x,C_r(\mathcal{X}_n\cap L_{\hat{f}_{h,n}}(\lambda)))<\eta$, which is a contradiction.

\end{document}